\documentclass[12pt]{amsart}%{article}

 \usepackage{amsfonts,amssymb,eucal}

\usepackage{amsthm}
 \usepackage{amsmath}
\usepackage{amscd}
 \usepackage{latexsym} 
\newcommand{\ac}{\ensuremath{\widetilde{{\text{SU}}(1,1)}}}
\newcommand{\hd}{\ensuremath{\mc{H}_{2,k}(\mc{D}_1)}}

 \newcommand{\e}{{\mbox{\rm e}}} 

  \newcommand{\mc}[1]{{\mathcal{#1}}}% mathematical caligraphic
 \newcommand{\got}[1]{{\mathfrak{#1}}}% gothic with mbox for  mathematic
\newcommand{\db}[1]{{\mathbb{#1}}}% double

\newcommand{\pa}{\partial}

\newcommand{\R}{\ensuremath{\mathbb{R}}}
\newcommand{\C}{\ensuremath{\mathbb{C}}}

 \newcommand{\Hi}{\ensuremath{\mathcal{H}}}% Hilbert space
 \newcommand{\Z}{\ensuremath{\mathbb{Z}}}

 \newtheorem{Remark}{Remark}

\newtheorem{Proposition}{Proposition}

\theoremstyle{definition}%nou 

\def\i{\mathrm{i}}
\begin{document}

\title{Generalized squeezed
   states for the Jacobi group}

\author{Stefan  Berceanu}
\address[Stefan  Berceanu]{National
 Institute for Physics and Nuclear Engineering\\
         Department of Theoretical Physics\\
         PO BOX MG-6, Bucharest-Magurele, Romania}
\email{Berceanu@theor1.theory.nipne.ro}

\begin{abstract}We analyze the relationship between the covering of
  the Jacobi group and the squeezed states.  We attach some
  nonclassical states to the Jacobi group.  The matrix elements of the
  Jacobi group are presented.
\end{abstract}
%\subjclass{81R30; 32Q15; 20C35; 81V80}
\keywords{Coherent states, squeezed states, Jacobi group}
\maketitle

%%%%%%%%%%%%%%%%%%%%%%%%%%%%%%%%%%%%%%%%%%%%
%% MAINMATTER
%%%%%%%%%%%%%%%%%%%%%%%%%%%%%%%%%%%%%%%%%%%%s

\section{Introduction\label{intr}}

The Jacobi group \cite{ez} -- the semidirect product of the
Heisenberg-Weyl group and the symplectic group -- is an important
object in the framework of Quantum Mechanics, Geometric Quantization,
Optics~\cite{stol1,stol2,stol3,stol4,ga1,ga2,ali,kbw2}. The Jacobi group
was investigated by physicists under other various names, as
``Schr\"odinger group'' \cite{nied} or ``Weyl-symplectic group''
\cite{kbw2}. The squeezed states \cite{stol1,stol2,stol3,stol4} in
Quantum Optics represent a physical realization of the coherent states
associated to the Jacobi group. Here we continue the investigation of
the Jacobi group started in \cite{jac1,sbj} using Perelomov coherent
states \cite{perG}.

  In \cite{jac1} we have constructed
generalized coherent states attached to the Jacobi group,
$G^J_1=H_1\rtimes \text{SU}(1,1)$, based on the homogeneous K\"ahler
manifold
$\mc{D}^J_1=H_1/\R\times\text{SU}(1,1)/\text{U}(1)=\C^1\times\mc{D}_1$.
Here $\mc{D}_1$ denotes the unit disk $\mc{D}_1=\{w\in\C||w|<1\}$, and
$H_n$ is the $(2n+1)$-dimensional real Heisenberg-Weyl (HW) group with
Lie algebra $\got{h}_n$. In \cite{jac1} we have also emphasized the
connection of our results with those of Berndt and Schmidt \cite{bs}
and K\"ahler \cite{cal}.  In \cite{sbj} we have considered coherent
states attached to the Jacobi group $G^J_n=H_n\rtimes
\text{Sp}(2n,\R)$, defined on the manifold
$\mc{D}^J_n=\C^n\times\mc{D}_n$, where $\mc{D}_n$ is the Siegel ball.

In the present note we underline the connection between the squeezed
states in Quantum Optics and the covering of the Jacobi group. 
 In
\cite{jac1} we have considered the action of the Jacobi group on the
minimal weight vector $e_0=\varphi_0\otimes\phi_{k0}$, where
$\varphi_0$ is the vacuum vector (${\bf a}\varphi_0=0$), while
$\phi_{k0}$ is the minimal weight vector of the positive discrete
series representation of the group $\text{SU}(1,1)$. The standard   
 squeezed
states correspond  to $e_0=\varphi_0$ \cite{stol1}.   In the
present paper we also give the matrix elements of the Jacobi group
acting on $\varphi_n\otimes\phi_{km}$ 
({\it generalized} coherent states attached to the
Jacobi group), where $\varphi_n$ is
$n$-particle (Fock) $\mc{F}$ vector, while $\phi_{km}$ is a normalized
vector obtained by the action of $(\bf{K}_+)^m$ on $\phi_{k0}$, $k\ge
0$, where $\bf{K}_+$ is the rising generator for the group
$\text{SU}(1,1)$.  More details about this calculation are given
elsewhere \cite{sac}. Many particular realizations
 of the generalized squeezed state for the Jacobi group are known.   We recall that the coherent states have been
introduced by Schr\"{o}dinger \cite{sch}, the squeezed states by
Kennard \cite{ken} and rediscovered later
\cite{stol1,stol2,stol3,stol4}, the displaced squeezed number states
by Husimi and Senitzky \cite{hs1,hs2}, the squeezed number states by
Plebanski \cite{pb}.

In the present note we follow the notation and convention of
\cite{jac1}. If $\pi$ is a representation of a Lie group $G$ with Lie
algebra $\got{g}$, then we denote ${\bf{X}}=d\pi (X),~ X\in\got{g}$.  

\section{The Jacobi group and its  algebra\label{unuu}}
% \subsection{The algebra}\label{jac1}

The Jacobi algebra is defined as the the semi-direct sum of the Lie
algebra $\got{h}_1$ of the Heisenberg-Weyl Lie group and the algebra
of the group ${\text{SU}}(1,1)$, $\got{g}^J_1= \got{h}_1\rtimes
\got{su}(1,1)$.

The Heisenberg-Weyl ideal $\got{h}_1= <\i s
1+xa^{\dagger}-\bar{x}a>_{s\in\R ,x\in\C}$ is determined by the
commutation relations
\begin{equation}
\left[a,{{K}}_+\right]=a^{\dagger}~;\left[ {{K}}_-,a^{\dagger}\right]=a;
\left[ {{K}}_+,a^{\dagger}\right]=\left[ {{K}}_-,a \right]= 0 ;
 \left[ {{K}}_0  ,a^{\dagger}\right]=\frac{1}{2}a^{\dagger};
 \left[ {{K}}_0,a\right]
=-\frac{1}{2}a,
\end{equation}
where ${a}^{\dagger}$ (${ a}$) are  the boson creation
(respectively, annihilation) operators,  which verify the canonical
commutation relations $[a,a^{\dagger}]=I, ~[a,I]= [a^{\dagger},I]=0$,
and $K_{0,+,-}$ are the generators of $\text{SU}(1,1)$ which satisfy
the commutation relations:
\begin{equation}
  \left[  {K}_{0},{K}_{\pm}\right]  =\pm {K}_{\pm}\ ,\quad\left[
    {K}_{-},{K}_{+}\right]
  =2{K}_{0}\ .\label{scs1}
\end{equation}

 We take a
 representation of $G^J_1$ (cf. \cite{jac1} and 
 (\ref{deplasare})-(\ref{gen1}) and  Proposition \ref{mm1} reproduced below)
 such that  the cyclic vector 
$e_0$  fulfills simultaneously the
conditions
\begin{equation}\label{cond}
{\bf{a}}e_0  =  0, 
~ {\bf{K}}_-e_0   =   0,~
{\bf{K}}_0e_0   =   k e_0;~ k>0,
\end{equation}
and we take $e_0=\varphi_0\otimes\phi_{k0}$. We consider for
$\text{Sp}(2,\mathbb{R})\approx \text{SU}(1,1)$ 
the  unitary irreducible positive
discrete series representation $D^+_k$ with Casimir operator $C={
  K}_{0}^{2}-{ K}_{1}^{2}-{ K}_{2}^{2}=k(k-1)$,
where $k$ is the Bargmann index for  $D^+_k$.
The orthonormal canonical basis of the representation space of
$\text{SU}(1,1)$ consists of the
vectors
\begin{equation}
\phi_{km} = \left[  \frac{\Gamma(2k)}{m!\Gamma(2k+m)}\right]  ^{1/2}\left(
{\bf{K}}_{+}\right)  ^{m}\phi_{k0}\ ,\label{scs2} m\in\Z_+ . %
\end{equation}
%\begin{equation}\label{iat}
%\bf{K}_{0}\phi_{k\,0}=k\phi_{k\,0}\ ,\quad \bf{K}_{-}\phi_{k\,0}=0\ .%
%\end{equation}

Also, in the Fock space $\mc{F}$,  we have 
$${\bf a}\varphi_n=\sqrt{n}\varphi_{n-1}; ~{\bf  a}^{\dagger}\varphi_n=
\sqrt{n+1}\varphi_{n+1}, ~\varphi_n=(n!)^{-\frac{1}{2}}({\bf
  a}^{\dagger})^{n}\varphi_0;~~
<\varphi_{n'},\varphi_n>=\delta_{nn'}.$$

Perelomov coherent state vectors associated to the group Jacobi
$G^J_1$, based on the manifold $\mc{D}^J_1$, are defined as
 \begin{equation}\label{csu}
e_{z,w}=\e^{z{\bf a}^{\dagger}+w{\bf{K}}_+}e_0, ~z,w\in\C,~ |w|<1 .
\end{equation}

We introduce the auxiliary operators
\begin{equation}\label{ssq}
{\bf{K}}_+  = \frac{1}{2}({\bf a}^{\dagger})^2+{\bf{K}}'_+ ,~
{\bf{K}}_-  = \frac{1}{2}{\bf a}^2+{\bf{K}}'_- ,~
{\bf{K}}_0  = \frac{1}{2}({\bf a}^{\dagger}{\bf a}+\frac{1}{2})+{\bf{K}}'_0 ,
\end{equation}
which have the properties
\begin{equation}\label{keg}
{\bf{K}}'_-e_0   =  0,~
{\bf{K}}'_0e_0 =  k'e_0; ~k=k'+\frac{1}{4}, 
\end{equation}
\begin{equation}
\left[{\bf{K}}'_{\sigma},{\bf a}\right]=\left[
{\bf{K}}'_{\sigma},{\bf a}^{\dagger}\right]=0,~\sigma =\pm ,0 ;~
 \left[ {\bf{K}}'_0,{\bf{K}}'_{\pm}\right]=\pm {\bf{K}}'_{\pm};~
 \left[ {\bf{K}}'_-,{\bf{K}}'_+\right]=2{\bf{K}}'_0 .
\end{equation}
The meaning of the splitting (\ref{ssq}) is explained
in  Theorem 2.6.1 from \cite{bs}, while the physical 
consequences  of this splitting are  briefly discussed  in \S
\ref{fizica}.
 More details are
given elsewhere \cite{sac}.  The positive discrete series corresponds
in (\ref{keg}) to $2k'=$integer.

We introduce the {\it displacement operator}
\begin{equation}\label{deplasare}
D(\alpha )=\exp (\alpha {\bf a}^{\dagger}-\bar{\alpha}{\bf a)}=
\exp(-\frac{1}{2}|\alpha
|^2)  \exp (\alpha {\bf a}^{\dagger})\exp(-\bar{\alpha}{\bf a}),
\end{equation}
and  the  {\it  unitary squeezed operator} of the
  $D^k_+$ representation of the group $\text{SU}(1,1)$,
 $\underline{S}(z)=S(w)$ ($w  =
\frac{z}{|z|}\tanh \,(|z|), \eta=\ln(1-|w|^2)$):
\begin{subequations}
\begin{eqnarray}
\underline{S}(z) & = &\exp (z{\bf{K}}_+-\bar{z}{\bf{K}}_-), ~z\in\C
;\label{u1} \\
S(w) & = &  \exp (w{\bf{K}}_+)\exp (\eta
{\bf{K}}_0)\exp(-\bar{w}{\bf{K}}_-), ~|w|<1 .\label{u2}
\end{eqnarray}
\end{subequations}
 We introduce also \cite{jac1}
the {\it generalized} squeezed coherent state vector
\begin{equation}\label{gen1}
\Psi_{\alpha, w}=D(\alpha )S(w) e_0,~
e_0=\varphi_0\otimes\phi_{k0}.
\end{equation}

We recall  some properties of the coherent states associated to the
group $G^J_1$, proved in \cite{jac1}:
  
\begin{Proposition}\label{mm1}
  The generalized squeezed coherent state
 vector {\em{(\ref{gen1})}}
and Perelomov coherent 
  state vector  {\em{(\ref{csu})}}
are related by the relation
\begin{equation}\label{csv}
  \Psi_{\alpha, w}= (1-w\bar{w})^k
  \exp (-\frac{\bar{\alpha}}{2}z)e_{z,w},\text{where}~
  z=\alpha-w\bar{\alpha}.
\end{equation}

Perelomov coherent state vector {\em{(\ref{csu})}} was calculated in
{\em{\cite{jac1}}} as
\begin{equation}\label{cemaie}
e_{z,w}=E(z,w)\varphi_0 \otimes\e^{w\bf{K}'_+}\phi_{k0},
\end{equation}
\begin{equation}\label{ezw}
  E(z,w)\varphi_0=\e^{z{\bf a}^{\dagger}+\frac{w}{2}({\bf a}^{\dagger})^2}\varphi_0
  =\sum_{n=0}^{\infty}\frac{P_n(z,w)}{(n!)^{1/2}}\varphi_n,
\end{equation}

\begin{equation}\label{marea}
  P_n(z,w)=n!\sum _{p=0}^{[\frac{n}{2}]}
  (\frac{w}{2})^p\frac{z^{n-2p}}{p!(n-2p)!} .
\end{equation}
The base of functions
$f_{nks}(\alpha,w)=<e_{\bar{\alpha},\bar{w}},\varphi_{n}\otimes\phi_{ks}>$,
where $k=k'+1/4,~2k'=$integer, $n,s=0,1,\cdots $, in which the
Bergman kernel {\em{(\ref{hot})}} can be expanded, is:
\begin{equation}\label{nlm}
  f_{nks}(\alpha,w)=f_{k's}(w)\frac{P_n(\alpha,w)}{\sqrt{n!}},f_{ls}(w)
  =\sqrt{\frac{\Gamma (s+2l)}{s!\Gamma (2l)}}w^s, ~ |w|<1.  
\end{equation}

The composition law in the Jacobi group $G^J_1=HW\rtimes SU(1,1)$ is
\begin{equation}\label{compositie}
(g_1,\alpha_1,t_1)\circ (g_2,\alpha_2, t_2)= (g_1\circ g_2,
g_2^{-1}\cdot\alpha_1+\alpha_2, t_1+ t_2 +\operatorname{Im}
(g^{-1}_2\cdot\alpha_1\bar{\alpha}_2)),
\end{equation}
where $g^{-1}\cdot\alpha =\bar{a}\alpha -b\bar{\alpha}$, and 
$g\in\text{SU}(1,1)$ is parametrized as 
\begin{equation}\label{ggg}
g=\left(\begin{array}{cc}a & b \\\bar{b} & \bar{a}\end{array}\right),
|a|^2-|b|^2=1, 
\end{equation}

Let $h = (g,\alpha )\in G^J_1$, $\pi (h)_k = T(g)_kD(\alpha )$, $g\in
SU(1,1)$, $\alpha\in \C$, and let $x=(z,w)\in
{\mc{D}^J_1}=\C\times\mc{D}_1$.  Then we have the formula:
\begin{equation}\label{rep}
\pi(h)_k\cdot   e_{z,w}= (\bar{a}+\bar{b}w)^{-2k}\exp(-\lambda_1)e_{z_1,w_1},
\end{equation}
\begin{equation}\label{x6}
\lambda_1=\frac{\bar{b}(z+z_0)^2}{2(\bar{a}+\bar{b}w)}
+\bar{\alpha}(z+\frac{z_0}{2}), ~ z_0= \alpha-\bar{\alpha} w,  
\end{equation}
\begin{equation}\label{xxx}
z_1=\frac{\alpha-\bar{\alpha}w+z}{\bar{b}w+\bar{a}}; ~ w_1
=\frac{a w+ b}{\bar{b}w+\bar{a}} . \end{equation}
The  space of functions $\Hi_K$ attached to the reproducing kernel
 \begin{equation}\label{hot}
K=
(1-w\bar{w})^{-2k}\exp{\frac{2z\bar{z}+z^2\bar{w}+\bar{z}^2w}{2(1-w\bar{w})}},  
 \end{equation}
consists of square integrable functions  with respect to the scalar
product 
\begin{equation}\label{ofi}
  (f,g)_k = \frac{4k-3}{2\pi^2} \! \int_{\alpha\in\C;
    |w|<1}\!\frac{\bar{f}(\alpha,w)g(\alpha,w)}{(1-|w|^2)^3K(w,\bar{w})}
  {\mathrm{d}}^{2}\alpha {\mathrm{d}}^{2}w.
\end{equation}
\end{Proposition}

\section{The covering of $\text{SU}(1,1)$\label{su11}}

If in (\ref{rep}) we take $\alpha= 0$, then
$T(g)_ke_w=(\bar{a}+\bar{b}w)^{-2k}e_{w_1}$, which corresponds to the
positive discrete series representation of $\text{SU}(1,1)$
\cite{bar47}
\begin{equation}\label{brrr}
[T(g)_k]f(z)=(\bar{a}+\bar{b}z)^{-2k}f(\frac{a z +b}{\bar{b}z+\bar{a}}),~
g=\left(\begin{array}{cc}a & b\\\bar{b} &\bar{a}\end{array}\right)\in\text{SU}(1,1),
\end{equation}
for the Hilbert space $\hd$, $2k=$ integer, 
 of holomorphic functions  on $\mc{D}_1$  with respect to the
scalar product \cite{bar47}
\begin{equation}\label{scc1}
(f,g )_{k}=\frac{2k-1}{\pi}
\int_{|w|<1} \bar{f}(w)g(w)(1-|w|^2)^{2k-2}{\mathrm{d}}^2w . 
\end{equation}

If in (\ref{rep}) we put $w=0$, we get  
$D(\alpha)e_z=\e^{-\frac{|\alpha|^2}{2}-z\bar{\alpha}}e_{z+\alpha}$,
which corresponds to  the Segal-Bargmann-Fock representation of the
Heisenberg group (see e.g. (1.71) in \cite{fol}).

 (\ref{rep}) corresponds to  the
continuous unitary representation $(\pi_k,\Hi_K)$ of $G^J_1$  
\begin{equation}\label{rep1}
(\pi(h)_k\cdot f)(x)=(\bar{a}+\bar{b}w)^{-2k}\exp(-\lambda_1)f(x_1),~  f\in\Hi_K,
\end{equation}
identified in \cite{jac1} with the result established in \cite{ez,bs}.

The functions $f_{nl}$ in (\ref{nlm}) form an orthonormal base in $\hd$,
$k>1/2$,  and if  
$f=\sum_{n=0}^{\infty} a_nz^n$, $g=\sum_{n=0}^{\infty} b_nz^n$, $|z|<1$, then 
\begin{equation}\label{newp}
(f,g)_k=\sum_{n=0}^{\infty}\frac{\Gamma (2k)\Gamma(n+1)}{\Gamma(2k+n)}\bar{a}_nb_n. 
\end{equation}

The representations of the group \ac~ have been considered in
\cite{puk,sal2_1,sal2_2,san}. For the universal covering group of
$\text{SU}(1,1)$, we use a class of Hilbert spaces indexed by a
parameter $k\in\R, ~0<|k|<1/2$ (cf. Sally \cite{sal2_1,sal2_2}).

Let now $\hd$ be the Hilbert space of functions
$f(z)=\sum_{n=0}^{\infty}a_nz^n$ holomorphic in $\mc{D}_1$ with norm
$$||f||_k^2=\sum_{n=0}^{\infty}\frac{\Gamma(2k)\Gamma(n+1)}{\Gamma(2k+n)}|a_n|^2<\infty
.$$

According to Bargmann ((4.3) in \cite{bar47}), starting from the matrix
 $g\in
\text{SU}(1,1)$ parame\-trized by (\ref{ggg}), an element $\tilde{g}$ of the covering group $\ac$ is parametrized by
$(\gamma,w)\in\R\times\mc{D}_1$ via  the relations
\begin{equation}\label{eqb}
a= e^{i\omega}(1-|\gamma|^2)^{-1/2}, ~ b
=e^{i\omega}\gamma (1-|\gamma|^2)^{-1/2}, |\gamma|<1,
\end{equation} 
and we have $\omega\in \left(-\pi/2,\pi/2\right]$, $\omega\in
\left(-\pi,\pi\right]$, $\omega\in \R$, for,
respectively, $\text{SO}^{\uparrow}(1,2)$,   
  $\text{Sp}(1,\R) \approx  \text{SU}(1,1)$, $\ac$. 

\begin{Remark} The representation {\em{(\ref{repac})}}
 of $\ac$ is  a continuous, irreducible  unitary representation  
with respect with the scalar product {\em{(\ref{newp})}}
\begin{equation}\label{repac}
[T(\tilde{g})_kf](z)=
e^{2ik\omega}(1-|\gamma|^2)^k(1+\bar{\gamma}z)^{-2k}f(\frac{a z
  +b}{\bar{b}z+
\bar{a}}),~ k\in\R, ~ 0<k<1/2.
\end{equation} 
  
\end{Remark} {\it Proof}. ~ Lemma (1.3.1) and Theorems (2.5.2),
(2.5.3) in the reference~\cite{sal2_1} are used. Here
we verify just that $(T(g)_kf_1,f_2)_k=(f_1,T(g^{-1})_kf_2)_k$, for
$g\in\text{SU}(1,1)$ and $f_{1,2}$ from the orthonormal set
(\ref{nlm}). We take $f_1=a_{kn}z^n$, $f_2=a_{kN}z^N$, where
$a^2_{kn}=\frac{\Gamma(n+2k)}{n!\Gamma(2k)}$.  Then $T(g)_k
f_1=A_{kn}(1-x)^{-q}(1+y)^n$, where
$A_{kn}=a_{kn}(\bar{a})^{-2k-n}b^n$, $x=-\lambda z, y=
\frac{z}{\lambda}$, $\lambda = \frac{b}{a}$, $q= 2k+n$, and
$g\in\text{SU}(1,1)$ is the matrix (\ref{ggg}). Then
$$(T(g)_kf_1,f_2)_k=\bar{A}_{kn}a_{kN}(\sum_{m=0}^{\infty}b_{qm}x^m\sum_{p=0}^nC^n_py^p,z^N)_k,$$
where  $b_{qm}=\frac{\Gamma(q+m)}{m!\Gamma(q)}$
and $C^n_k=\frac{n!}{k!(n-k)!}$. 
With the orthogonality of the system (\ref{nlm}), we get
 $$(T(g)_kf_1,f_2)_k=\frac{\bar{A}_{kn}}{a_{kN}}\sum_{m=0}^{\infty}(-\lambda)^mb_{qm}\sum_{p=0}^n
C^n_p\bar{\lambda}^{-p}\delta_{m+p,N} . $$
Similarly, we have ($Q=2k+N$)
$$(f_1,T(g^{-1})_kf_2)_k=
\frac{a_{kN}}{a_{kn}}a^{-Q}(-b)^N\sum_{M=0}^{\infty}b_{QM}(\frac{\bar{b}}{a})^M\sum_{p=0}^NC^N_p
(-\frac{\bar{a}}{b})^p\delta_{n,p+M}.$$
So, we have (we take $N-n\ge 0$)
 $$(T(g)_kf_1,f_2)_k=\frac{a_{kn}}{a_{kN}}a^{-q}\bar{a}^N\bar{b}^{n-N}\sum_{m=N-n}^N b_{qm}C^{n}_{N-m}
(-|b/a|)^{2m},$$
$$(f_1,T(g^{-1})_kf_2)=\frac{a_{kN}}{a_{kn}}a^{-Q}(-b)^N(-\frac{\bar{a}}{b})^n\sum_{M=0}^n
b_{QM}C^N_{n-M}
(-|b/a|)^{2M}.$$
With the change of variables $m-(N-n)=M$, it is easy to check that the
last two expressions are identical.

\begin{Remark}\label{DIR}It can be checked out  that the differential operators
 $\db{K}_-=\frac{\pa}{\pa z},~\db{K}_0=k+
z\frac{\pa}{\pa z},~
\db{K}_+=2kz +z^2\frac{\pa}{\pa z}$ 
 corresponding to the generators of $\text{SU}(1,1)$ 
have the adequate  hermitian conjugate properties
 with respect to the scalar product
$(\cdot,\cdot)_k$ \emph{(\ref{newp})}, $k>0$. 
\end{Remark}
\section{Squeezed states and Jacobi group\label{fizica}}

The standard squeezed states \cite{stol1,stol2,stol3,stol4} correspond 
to the action of Jacobi group on the extremal weight vector 
$e_0=\varphi_0$.  This corresponds to take zero the
action of ${\bf K}'$ in the splitting (\ref{ssq}). 
 This is the group which admits the so
called {\it Schr\"odinger-Weil} representation $\pi^m_{SW}$ with
character $\psi^{m}(x)=\e^{2\pi\i mx}, ~m\in\R$, considered in the
mentioned Theorem 2.6.1 in \cite{bs}, where $x$ is in the center of
the Heisenberg-Weyl algebra. The part of the representation
corresponding to the covering group $\tilde{G}^J_1$ is called ``Weil
representation'' at p. 23 in \cite{bs}.

We see that we have
${\bf{K}}_0\varphi_{2p}=(p+\frac{1}{4})\varphi_{2p},
 ~
{\bf{K}}_0\varphi_{2p+1}=(p+\frac{3}{4})\varphi_{2p+1}$, and 
irreducible representations with $k=\frac{1}{4}$,
$k=\frac{3}{4}$ of {\it  the covering group  \ac~ must be considered}.
 The vacuum squeezed
state contains only even Fock states.

So, dealing with squeezed states,  we have  to consider   the covering of the  Jacobi group $G^J_1$, 
$ \tilde{G}^J_1=HW\rtimes \ac $.

 {\it 
The orthonormal basis of $\Hi_K$ for $G^{J}_1$ in the realization
${\bf{K}}_+
=\frac{1}{2}({\bf{a}}^{\dagger})^2$ in the splitting {\em{(\ref{ssq})}}, 
 where the reproducing
kernel $K$ is given by {\em{(\ref{hot})}}
with $k=1/4$,  consists of polynomials
$(n!)^{-1/2}P_n$, $n=0,1,\cdots$. Instead of the formula
{\em{(\ref{rep})}} for  $G^J_1$ ($2k'$ = integer), we get  
for   $\tilde{G}^J_1$ ($k'>0$) the formula}
\begin{equation}\label{fanta}
\pi(\tilde{h})_k\cdot e_{z,w}=
\e^{2\i k\omega}(1-|\gamma|^2)^k(1+\bar{\gamma}w)^{-2k}\exp(- \lambda_1)
e_{z_1,w_1}. 
\end{equation}

%%%%%\end{Proposition}

\section{Matrix elements for the  Jacobi group}
%\subsection{Notation}

Now we briefly present the main steps in the calculation of the matrix
elements of $G^{J}_1$ with respect to the considered representation:
\begin{equation}\label{MNM}
\left\langle \varphi_{n^{\prime}}\otimes\phi_{km^{\prime}}\right\vert
D(\alpha)S\left(  w\right)  \left\vert \varphi_{n}
\otimes\phi_{km}\right\rangle .\end{equation}
More details are  given elsewhere \cite{sac}.

Firstly, we re-obtain the matrix elements of the Heisenberg-Weyl group \cite{fey,sw}: 
\begin{equation}
\left\langle \varphi_{m}\right\vert D(\alpha)\left\vert \varphi_{n}%
\right\rangle   =\sqrt{\frac{n!}{m!}}\alpha^{m-n}L_{n}^{m-n}(|\alpha
|^{2})\exp\left(  -\left\vert \alpha\right\vert ^{2}/2\right)  \ ,\quad m\geq
n\ ,\label{s6}
\end{equation}
where  $m$,  $n$ are  non-negative integers and  $L_{n}^{s}$ are the    associated Laguerre
polynomials.
%\begin{equation}
%L_{n}^{s}(x)=\sum_{k=0}^{n}\frac{(n+s)!(-x)^{k}}{k!(n-k)!(k+s)!}\ .\label{s7}%
%\end{equation}

%{\bf  The hypergeometrical function}
%{\bf Bargmann's matrix elements}
 Next we calculate the matrix elements
 $S(w)_{km'm} =\left\langle \phi_{k\,m^{\prime}}\right\vert S\left(  w\right)  \left\vert
\phi_{k\,m}\right\rangle $. We use the relations:
\begin{subequations}\label{b11}
\begin{eqnarray}
{\bf{K}}_{0}\phi_{k\,m}  & = & (k+m)\phi_{k\,m},\label{oo1} \\
{\bf{K}}_{+}\phi_{k\,m}  & = & \left[  (m+1)(m+2k)\right]  ^{1/2}\phi_{k\,m+1}
\ ,\label{oo2}\\
{\bf{K}}_{-}\phi_{k\,m}  & = & \left[  m(m+2k-1)\right]  ^{1/2}\phi_{k\,m-1}\ ,\quad
m>0, m\in\Z_+ . \label{oo3} 
\end{eqnarray}
\end{subequations}
We introduce (\ref{oo3}) in the expression of $S(w)$ and we get
\begin{align}
\exp\left(  -\bar{w}{\bf K}_{-}\right)  \phi_{k\,m} &  =%
%TCIMACRO{\dsum \limits_{p=0}^{m}}%
%BeginExpansion
{\displaystyle\sum\limits_{p=0}^{m}}
%EndExpansion
\frac{\left(  -\bar{w}\right)  ^{p}}{p!}{\bf K}_{-}^{p}\phi_{k\,m}\label{h3}\\
&  =%
%TCIMACRO{\dsum \limits_{p=0}^{m}}%
%BeginExpansion
{\displaystyle\sum\limits_{p=0}^{m}}
%EndExpansion
\left[  \frac{m!}{(m-p)!}\frac{\Gamma(2k+m)}{\Gamma(2k+m-p)}\right]
^{1/2}\frac{\left(  -\bar{w}\right)  ^{p}}{p!}\phi_{k\,m-p} .\nonumber
\end{align}%

Then we apply successively (\ref{oo1}) and (\ref{oo2}) and we have ( $m^{\prime}\geq m$):
\begin{eqnarray}\label{h11}
S(w)_{km'm} &   = &  
\left[  \frac{m^{\prime}!\Gamma(2k+m^{\prime})}{m!\Gamma(2k+m)}\right]
^{1/2}
\frac{w^{m^{\prime}-m}}{(m^{\prime
}-m)!}\left(  1-w\bar{w}\right)  ^{k+m}\\
~~ & \times & 
F\left(  -m,1-2k-m;m^{\prime}-m+1;\frac{-w\bar{w}}{1-w\bar{w}}\right).
\nonumber%
\end{eqnarray}

The matrix elements of the group    $\text{SU}(1,1)$  for the discrete
series representations have been calculated in \cite{bar47},
and the case corresponding to the \ac~ for $k=1/4, 3/4$ in \cite{pp}.

Using Kummer's formula
$F\left(  a,b;c;\frac{z}{z-1}\right)  =\left(  1-z\right)  ^{a}
F(a,c-b;c;z)$, we put (\ref{h11}) in the form
 \begin{equation}\label{adoua}
   S(w)_{km'm}   \! =\! 
   \left[  \frac{m^{\prime}!\Gamma(2k+m^{\prime})}{m!\Gamma(2k+m)}\right]
   ^{\frac{1}{2}}\!\frac{w^{m^{\prime}-m}\left( \! 1\!-\!w\bar{w}\right)  ^{k}}{(m^{\prime
     }-m)!}
   F\left( \! -m,\!2k\!+\!m^{\prime};m^{\prime}\!-\!m\!+\!1\!;\!|w|^2\right)
\end{equation}
%\begin{eqnarray}\label{h17}
%S(w)_{km'm} &   = &  
%\left[  \frac{m^{\prime}!\Gamma(2k+m^{\prime})}{m!\Gamma(2k+m)}\right]
%^{1/2}\frac{w^{m^{\prime}-m}\left(  1\!-\!w\bar{w}\right)  ^{k+m}}{(m^{\prime
%}-m)!}\\
%~~ & \times & 
%F\left(  -m,2k+m^{\prime};m^{\prime}\!-\!m\!+\!1;|w|^2\right) ,
%\nonumber%
%\end{eqnarray} 
for $m^{\prime}\geq m$, while for 
  $m\geq m^{\prime}$, we have 
$S_{km^{\prime}m}(w)=S_{kmm^{\prime}}(-\bar{w})\label{h15}$. 

Now we calculate the  matrix elements (\ref{MNM}). Starting
from the splitting (\ref{ssq}),  we
introduce the auxiliary operator
%\bigskip%
\begin{align}
S\left(  w,w^{\prime}\right)    & =\exp\left[ \frac{w}{2}(a^{\dagger})^{2}\right]
\exp\left[  \frac{\eta}{2}(a^{\dagger}a+\frac{1}{2})\right]  \exp\left(
-\frac{\bar{w}}{2}a^{2}\right) \\ \nonumber & \otimes\exp\left(  w^{\prime}K_{+}^{\prime
}\right)  \exp\left(  \eta^{\prime}K_{0}^{\prime}\right)  \exp\left(
-\overline{w^{\prime}}K_{-}^{\prime}\right), \\
\eta & =\ln\left(  1-w\bar{w}\right)  ,\quad\eta^{\prime}=\ln\left(
1-w^{\prime}\overline{w^{\prime}}\right) . \nonumber
\end{align}

{\it The matrix elements $S_{km^{\prime}m}(w)$ are given
  by {\em{(\ref{h11})}}, {\em{(\ref{adoua})}}. {\em{(\ref{h11})}} is a particular case of
    formula {\em{(10.28)}} in {\em{\cite{bar47}}}.  Moreover, the formula 
{\em{(10.28)}}  in {\em{\cite{bar47}}}  can be re-obtained  from our 
    {\em{(\ref{h11})}}  using the relation 
\begin{equation}\label{BBB}
T (g) _k\phi_{km}=
\left(\frac{a}{|a|}\right)^{2(k+m)}
S\left(\frac{b}{\bar{a}}\right)\phi_{km},~g=\left( \begin{array}{cc} a & b\\
  \bar{b} &\bar{a} \end{array}\right)\in{\text{SU}}(1,1) .
\end{equation}

  Using the identification  $\varphi_{2n+\epsilon}\leftrightarrow\phi_{1/4+\epsilon/2;n}$,
($\epsilon=0$ or 1), the    matrix elements {\em{(\ref{MNM})}} of the
Jacobi group $G^{J}_1$ are
obtained taking $w=w'$ in formula {\em{(\ref{n15})}}: 
\begin{align}
\label{n15} \left\langle \varphi_{n^{\prime}}\otimes\phi_{km^{\prime}}\right\vert
D(\alpha)S\left(  w,w^{\prime}\right)  \left\vert \varphi_{2s+\varepsilon
}\otimes\phi_{km}\right\rangle  &  = S_{km^{\prime}m}(w^{\prime})\\%
%TCIMACRO{\dsum \limits_{s^{\prime}\geq0}}%
%BeginExpansion
\nonumber & \times {\displaystyle\sum\limits_{s^{\prime}\geq0}}
%EndExpansion
 \left\langle \varphi_{n^{\prime}}\right\vert D(\alpha)\left\vert
\varphi_{2s^{\prime}+\varepsilon}\right\rangle S_{1/4+\varepsilon
/2\ s^{\prime}s}(w).
\end{align}%
 The matrix elements $\pi(h)$ are obtained from {\em{(\ref{n15})}}
  taking into
account {\em{(\ref{BBB})}}.}

\end{document}